\documentclass[a4paper,10pt]{article}
\usepackage{fullpage}
\usepackage{amsmath,amssymb,latexsym,amsthm}
\newcommand{\R}{\mathbb{R}}

\renewcommand{\geq}{\geqslant}
\renewcommand{\div}{\operatorname{div\,}}
\newcommand{\curl}{\operatorname{curl}}

\newcommand{\Id}{\operatorname{Id}}
\def\eop{\hfill $\Box$ \\ \ \par}
\newtheorem{Theorem}{Theorem}
\newtheorem{Definition}{Definition}

\newtheorem{Proposition}{Proposition}
\newtheorem{Lemma}{Lemma}
\newtheorem{Remark}{Remark}

\title{The movement of a solid in an incompressible perfect fluid \\ as a geodesic flow}
\author{Olivier Glass\footnote{CEREMADE,
Universit\'e Paris-Dauphine, 
Place du Mar\'echal de Lattre de Tassigny,
75775 Paris Cedex 16,
FRANCE
},
Franck Sueur\footnote{Laboratoire Jacques-Louis Lions,
Universit\'e Pierre et Marie Curie - Paris 6,
4 Place Jussieu,
75005 Paris,
FRANCE
}, }
\begin{document}
\maketitle
\begin{abstract}
The motion of a rigid body immersed in an incompressible perfect fluid which occupies a three-dimensional bounded domain have been recently studied under its PDE formulation. In particular  classical solutions have been shown to exist locally in time.
In this note, following the celebrated result of Arnold \cite{Arnold} concerning the case of a perfect incompressible fluid alone, we prove that these classical solutions are the geodesics of a Riemannian manifold of infinite dimension, in the sense that they are the critical points of an action, which is the integral over time of the total kinetic energy of the fluid-rigid body system. 
\end{abstract} \par
\ \par
\noindent
{\small {\bf Keywords.} Perfect incompressible fluid, fluid-rigid body interaction, least action principle.\\
{\bf AMS Subject Classification. } 76B99, 74F10.} \\
%
%
%
%
%%%%%%%%%%%%%%%%%%%%%%%%%%%%%%%%%%%%%%%%%%%%%%%%%%%%%%
%
% 

\section{Introduction}
\label{Intro}

We consider  the motion of a rigid body immersed in an incompressible homogeneous perfect fluid, so that the system fluid-rigid body occupies a smooth open and bounded domain $\Omega \subset \R^{3}$. The solid is supposed to occupy at each instant $t \geq 0$ a smooth closed connected subset $\mathcal{S}(t) \subset \Omega$ which is surrounded by a perfect incompressible fluid filling   the domain $\mathcal{F}(t) :=   \Omega \setminus  \mathcal{S}(t)$. \par
\ \par
For the point of view of PDEs, this system have been recently studied in \cite{ort1}, \cite{ort2}, \cite{rosier}, \cite{ht}, \cite{ogfstt} which have set a Cauchy theory for classical solutions. \par
\ \par
The aim of this note is to provide a rigorous proof that the classical solutions can be equivalently thought as  geodesics of a Riemannian manifold of infinite dimension, in the sense that they are the critical points of an action, which is the integral over time of the total kinetic energy of the fluid-rigid body system. It was pointed out in a famous paper by Arnold \cite{Arnold} that both the Euler equations for a rigid body as well as the Euler equations for a perfect fluid can be derived with this approach. 
The motion of a rigid body in a frame attached to its center of mass can be considered as a geodesic on the special orthogonal group $SO(3)$. On the other hand the motion of a perfect fluid filling a container $ \Omega$ (without any immersed rigid body in it) can  be considered as a geodesic equation on the space $\text{Sdiff}^{+} (\Omega)$ of the volume and orientation preserving diffeormorphisms of $\Omega$. \par
\ \par
It is hence natural to try to extend this analysis to a system of interaction of a perfect fluid and a rigid body. In particular cases, when the fluid is irrotational or when the vorticity of the fluid is given by a finite number of point vortices, so that the dynamics is finite-dimensional, this was studied in details in \cite{VKM1,VKM2}; see also references therein. The goal of this paper is to observe that one can see the motion of a rigid body in a fluid governed by the incompressible Euler equations, as a geodesic flow, in the presence of a regular distributed vorticity, as well. \par
\ \par
The structure of this paper is as follows. In Subsection \ref{Pde}, we first recall the PDE formulation of the system. Then in Subsection \ref{GEO}, we describe the infinite-dimensional manifold and the action used in the geometric formulation of the problem. In Subsection \ref{Equi}, we state the main result of the paper, that is, the equivalence of the two points of view. Section \ref{Sec:Proofs} contains the proofs of the various statements. \par
\ \par
\subsection{PDE formulation} 
\label{Pde}
 The dynamics of this system can be described thanks to the following PDEs:
\begin{eqnarray}
\label{Euler1a2}
\displaystyle \frac{\partial u}{\partial t}+(u\cdot\nabla)u + \nabla p &=& 0  \ \text{for}  \ x\in   \mathcal{F}(t) ,  \\
\label{Euler2a2}
\div u &=& 0   \ \text{for}  \ x\in   \mathcal{F}(t), \\ 
\label{Solide1} 
m x_{B}''(t) &=&  \int_{\partial \mathcal{S}(t)} pn \ d\Gamma, \\
\label{Solide2}
(\mathcal{J} r )'(t) &=&  \int_{\partial \mathcal{S}(t)} (x- x_{B})\times pn \ d\Gamma  , \\
\label{Euler3a2}
u\cdot n &=& 0 , \ \text{for}  \ x\in \partial \Omega ,  \\
\label{Euler3b} 
u\cdot n &=& v \cdot n , \ \text{for}  \    \ x\in \partial \mathcal{S}(t),  \\
\label{Eulerci2}
u |_{t= 0} &=& u_0 , \\
\label{Solideci1}
\ell(0) &=& \ell_0, \\
\label{Solideci2}
r  (0)&=& r _0, \\
\label{CMci}
x_{B}(0) &=& x_0 .
\end{eqnarray}
Equations (\ref{Euler1a2})  and (\ref{Euler2a2}) are the incompressible Euler equations.
The vector field $u$ is the fluid velocity and the scalar field $p$ denotes the pressure.
Equations  (\ref{Solide1})  and (\ref{Solide2}) are Newton's laws for linear and angular momenta of the body under the influence of the pressure force. Here we denote by $m$ the mass of the rigid body (normalized in order that the density of the fluid is $\rho_{F} = 1$), by $x_B (t)$  the position of its center of mass, 
$n(t,x)$ denotes the unit normal vector pointing outside the fluid
and $ d\Gamma(t) $ denotes the surface measure on $\partial \mathcal{S}(t)$. The time-dependent vector 
\begin{equation} \label{DefEll}
\ell(t) :=  x_{B}'(t),
\end{equation}
denotes the velocity of the center of mass of the solid and $r$ denotes its angular speed, so that the solid velocity field is given by
\begin{equation} \label{vietendue}
v(t,x) :=\ell(t) +   r (t) \times (x -  x_{B}(t)).
\end{equation}
In (\ref{Solide2}) the $3 \times 3$ matrix $\mathcal{J}={\mathcal J}(t)$ denotes the moment of inertia which depends on time according to Sylvester's law
\begin{equation}\label{Sylvester}
 \mathcal{J} =Q \mathcal{J}_0 Q^{*} ,
 \end{equation}
where $\mathcal{J}_0$ is the initial value of $\mathcal{J}$ and  where the rotation matrix $Q \in SO(3)$ is deduced from $r$ by the following differential equation (where we use the convention to consider the operator $r(t) \times \cdot$ as a matrix):
\begin{equation} \label{LoiDeQ}
Q'(t) = r(t) \times Q(t) \text{ and } Q(0) = \Id_3.
\end{equation}
The matrix ${\mathcal J}_{0}$ can be obtained as follows. Given a positive function $\rho_{{\mathcal S}_0} \in L^{\infty}({\mathcal S}_{0};\R)$ describing the density in the solid (again normalized in order that the density of the fluid is $\rho_{F} = 1$), the data $m$, $x_{0}$ and ${\mathcal J}_{0}$ can be computed by it first moments
\begin{equation} \label{EqMasse}
m :=  \int_{{\mathcal S}_0} \rho_{{\mathcal S}_0} dx  > 0,
\end{equation}
\begin{equation} \label{Eq:CG}
m x_0 :=   \int_{{\mathcal S}_0} x \rho_{{\mathcal S}_0} (x) dx,
\end{equation}
\begin{equation}\label{eqJ}
\mathcal{J}_{0} := \int_{ \mathcal{S}_{0}}  \rho_{{\mathcal S}_0}(x)   \big( | x- x_{0} |^2 \Id_3 -(x- x_{0})  \otimes   (x- x_{0})    \big) dx  .
\end{equation}
Finally, the domains occupied by the solid and the fluid are given by
\begin{equation} \label{PlaceDuSolideEtDuFluide}
{\mathcal S}(t) = \Big\{ x_{B}  (t) + Q(t)(x-x_{0}), \ x \in {\mathcal S}_{0} \Big\} \text{ and } {\mathcal F}(t)=\Omega \setminus \overline{{\mathcal S}(t)},
\end{equation}
starting from a given initial position  $\mathcal{S}_0 \subset \Omega$, such that $\mathcal{F}_0 :=  \Omega \setminus  \mathcal{S}_0$. Let us underline that ${\mathcal S}(t)$ and $\partial \Omega$ being compact, and since ${\mathcal S}(t) \subset \Omega$, the solutions that we consider satisfy
\begin{equation} \label{DistanceAuBord}
\mbox{dist}({\mathcal S}(t),\partial \Omega) >0 \text{ on } [0,T].
\end{equation}
Let us give a precise definition of the classical solutions examined in this paper.
\begin{Definition}[Classical solutions] 
 \label{germi}
 We call  classical solution of the PDEs formulation on $[{0}, T]$ some
\begin{equation*}
(u,x_B , r)   \in  C^{1,\lambda} (\cup_{t \in [0,T]} \Big( \{t\} \times {\mathcal F}(t) \Big) ;\R^{3})  \times  C^{2} ([0,T];\R^{3}) \times C^{1} ([0,T];\R^{3}),
\end{equation*}
(for some $\lambda \in (0,1)$) satisfying \eqref{Euler1a2}--\eqref{DistanceAuBord}.
\end{Definition}
The local-in-time existence and uniqueness of classical solutions to the problem \eqref{Euler1a2}--\eqref{DistanceAuBord} holds when the initial velocity of the fluid  is in the H\"older space $C^{1,r}$ cf. \cite{ogfstt} for a precise statement. Let us also mention the earlier results of  Ortega, Rosier and Takahashi \cite{ort1}-\cite{ort2} where the body-fluid system occupies the plane $\R^{2}$, Rosier and Rosier \cite{rosier} in the case of a body in $\R^{3}$ and  Houot, San Martin and Tucsnak \cite{ht} in the case (considered here) of a bounded domain, with the  initial velocity in a Sobolev space $H^{m}$, $m \geq 3$.

\subsection{Geodesic formulation} 
\label{GEO}
Let us now turn to the geometric viewpoint. We first describe below the infinite-dimensional space of configuration of the system. Next we introduce a natural action, which allows to define our notion of geodesic.
\subsubsection{Rigid movements}
Let us first describe the rigid part of the motion. To a velocity vector field $v$ of a rigid body one associates the flow 
\begin{equation}
\label{flowS}
\partial_t \tau(t,x)  = v (t,  \tau(t,x) ) \text{ and } \tau(0,x) = x \text{ for } (t,x) \in [0,T] \times {\mathcal S}_{0}.
\end{equation}
It is easy to integrate to find 
\begin{equation*}
\tau(t,x) = x_{B}  (t) + Q(t)(x-x_{0}),
\end{equation*}
where 
\begin{equation} \nonumber
x_{B}  (t) := \tau(t,x_0 ),
\end{equation}
and where $ Q(t)$ is obtained from $v$ by \eqref{LoiDeQ} and $r$ is given by \eqref{vietendue}. \par
The flow $\tau$ can be seen as a $C^1$ function of the time with values in the Lie group $SE(3)$ of rigid motions (the special Euclidean group), that is the group generated by translations and rotations in $3$D. Its tangent space at $\Id \in SE(3)$ is
\begin{eqnarray*}
\mathfrak{se}(3) := T_{\Id} SE(3) =  \Big\{ v \in C^{1}(\R^{3};\R^{3})  \ \Big/ \ D(v) = 0 \Big\}, 
\end{eqnarray*}
where  $D(v)$ denotes the tensor of deformations $2 D(v) := (\partial_i v_j +  \partial_j v_i )_{i,j} $.
Given $x_{B}$ in $\R^3$, we have the following:
\begin{equation*}
\mathfrak{se}(3)  = \Big\{ v :\R^{3} \rightarrow \R^{3}  \ \Big/ \  \exists (\ell,r) \in \R^{3} \times  \R^{3}  , \ \forall x \in \R^3 , \
v(x) = \ell + r \times ( x- x_B ) \} .
\end{equation*}
Moreover given $x_{B}$ in $\R^3$ and $v \in \mathfrak{se}(3) $,  the ordered pair $(\ell,r)$ above is unique. Hence $\mathfrak{se}(3) \sim \R^{3} \times \mbox{Skew}_{3}(\R)$. \par 
Accordingly, the tangent space of $SE(3)$ at $\tau \in SE(3)$ is 
\begin{equation*}
T_\tau SE(3) :=  \Big\{ v  \circ \tau  \text{ with } v \in \mathfrak{se}(3) \Big\}.
\end{equation*}
Now, $x_0$ being given,  we introduce the following projections on $T_{\tau} SE(3)$: to $\sigma = v \circ \tau \in T_{\tau} SE(3)$ we associate 
$(L_{\tau}[\sigma] , R_{\tau}[\sigma])$ the unique ordered pair $(\ell,r)$ associated to $v$ with $x_{B} := {\tau} (x_0)$, in other words:
\begin{equation*}
\sigma(\tau^{-1}(x)) = v(x)= L_{\tau}[\sigma] + R_{\tau}[\sigma] \times (x - \tau(x_{0})).
\end{equation*}
\ \par
Let us conclude this subsection by describing the energy of the solid. Using the choice of $x_{B} (t)$ as the center of mass of the body at time $t$, we have that 
\begin{equation}
\label{premeti}
 \int_{{\mathcal S}(t) }  \rho_{{\mathcal S}_0} (\tau_t^{-1}(x)) \, (x- x_{B}(t)) ) \, dx = 0,
\end{equation}
 and therefore,   for any $(\ell_1 , r_1 ,\ell_2 , r_2 ) \in (\R^3)^4$, for any $t$, 
\begin{equation}\label{meti}
 \int_{{\mathcal S}(t) }  \rho_{{\mathcal S}_0}(\tau_t^{-1}(x)) \, (\ell_1 +   r_1  \times (x -  x_{B}(t)) ) \cdot (\ell_2 +   r_2  \times (x -  x_{B} (t)) ) \, dx =  m \ell_1 \cdot \ell_2 +  \mathcal{J} (t) r_1 \cdot r_2 ,
  \end{equation}
where $ \mathcal{J} (t) $ is given by \eqref{Sylvester} and 
the notation $\tau_t^{-1}$ stands for the  inverse of the function $\tau_{t}:=\tau (t,\cdot)$.
\subsubsection{Fluid displacements and Arnold's geodesic interpretation} 
\label{Arnold}
Let us briefly recall Arnold's interpretation of the Euler equation.
To a velocity vector field $u$ satisfying the incompressible Euler equations in $\Omega$ (without body) one associates the flow $ \eta$ defined on $[0,T] \times \Omega$ by 
\begin{equation}
\label{flowF}
\partial_t  \eta(t,x)  = u (t, \eta(t,x) )  \text{ and }
\eta(0,x) = x .
\end{equation}
The flow $\eta$ can be seen as a continuous function of the time with values in the  space $\text{Sdiff}^{+}  (\Omega)$ of the volume and orientation preserving diffeormorphisms defined of $\Omega$. The latter is viewed as an infinite-dimensional manifold with the metric inherited from the embedding in $L^2  (\Omega ; \R^{3})$, and the tangent space in $\eta \in \text{Sdiff}^{+}  (\Omega)$ is 
\begin{equation*}
T_\eta \ \text{Sdiff}^{+}  (\Omega) :=  \Big\{ u \circ \eta  \text{ with } u  \in  C^{1} ( \Omega  ;\R^{3})  \text{ such that }  \div(u)=0 \ \text{ in }\ \Omega \ \text{ and }\ u \cdot n=0 \ \text{ on }\ \partial \Omega  \Big\}.
\end{equation*}
Euler's equations are then interpreted as a geodesic equation on $\text{Sdiff}^{+}  (\Omega)$.
The pressure field appears as a Lagrange multiplier for the divergence-free constraint on the velocity.
Ebin and Marsden proved in \cite{EbinMarsden} the existence of these geodesics when the initial velocity of the fluid  is in the Sobolev space $H^{s} (\Omega)$, $s>\frac{5}{2}$. \par
\subsubsection{Possible configurations of the fluid/body system as a Riemannian manifold}
\label{PossibleConfigurations}
To introduce the geodesic formulation of the motion of a rigid body immersed in an incompressible perfect 
fluid, we first describe the infinite-dimensional manifold on which these geodesics will be considered. \par
\ \\
{\bf The set of the possible configurations as a Riemannian manifold.} To begin with, we first describe the set of the possible configurations of the system at a fixed time by setting
\begin{multline} \nonumber
{\mathcal C} := \Big\{ \ (\tau,\eta) \in SE(3) \times C^{1,\lambda}( {\mathcal F}_{0}; \R^{3})
\text{ such that } \tau({\mathcal S}_{0}) \subset \Omega, \\
\eta \text{ is a  volume and orientation preserving diffeomorphism } {\mathcal F}_{0} \rightarrow \Omega \setminus \left[\tau({\mathcal S}_{0})\right] \Big\}.
\label{defC}
\end{multline}
We will represent $(\tau,\eta)$  by $\phi : \Omega \rightarrow \Omega$ such that  $\phi_{|{\mathcal S}_{0}}= \tau$ and $\phi_{|{\mathcal F}_{0}}=\eta$. Note that $\phi$ is not necessarily continuous. \par
Let us observe that according to the Helmoltz decomposition, in $C^{1,\lambda}( {\mathcal F}_{0}; \R^{3})$, the space of divergence-free vector fields is closed and admits a topological complement. Therefore, one can see that ${\mathcal C}$ is a submanifold of a manifold modelled on the Banach space
\begin{equation*} 
E :=  \mathfrak{se}(3)  \times C^{1,\lambda}( {\mathcal F}_{0}; \R^{3}).
\end{equation*}
We will be interested in the tangent space of this infinite-dimensional manifold. Let us first recall that by definition the tangent space $T_{(\tau,\eta)} {\mathcal C}$ in $(\tau,\eta) \in {\mathcal C}$ is the set of equivalence classes of germs of  continuously differentiable function $\Theta: (-\varepsilon, \varepsilon) \rightarrow {\mathcal C}$ (for some $\varepsilon>0$) such that $\Theta(0)=(\tau,\eta)$, for the equivalence relation:
\begin{equation*}
\Theta_{1} \sim \Theta_{2} \text{ iff } \Theta_{1}'(0) =\Theta_{2}'(0)  \text{ in local charts}.
\end{equation*}
Equivalently $(\sigma,\mu) \in T_{(\tau,\eta)} {\mathcal C}$ is the set of equivalence classes of pairs $(\psi, e)$ where $\psi$ is a chart defined on a neighborhood of $(\tau,\eta) \in {\mathcal C}$ and $e \in E$ for the relation
\begin{equation*}
(\psi_1 ,e_1 ) \equiv (\psi_2 ,e_2 ) \ \text{ if } \
 D ( \psi_2 \circ \psi_1^{-1} )(  \psi_1 (\tau,\eta ) ) \cdot e_1 = e_2 .
\end{equation*}
Clearly this makes  $T_{(\tau,\eta)} {\mathcal C}$  a linear space and  $T {\mathcal C} := \displaystyle \bigcup_{(\tau,\eta)\in  {\mathcal C}} \Big( \{(\tau,\eta)\} \times T_{(\tau,\eta)} {\mathcal C} \Big)$ a vector bundle whose base space is ${\mathcal C}$. \par
\ \par
\noindent
Let us introduce the following notation: given $(\tau,\eta) \in {\mathcal C}$ and $\mu \in C^{1,\lambda}( {\mathcal F}_{0}; \R^{3})$, we introduce $U_{\eta}[\mu]: \eta({\mathcal F}_{0}) \rightarrow \R^{3}$ by
\begin{equation} \nonumber
U_{\eta}[\mu] := \mu \circ \eta^{-1}.
\end{equation}
Recall that due to the openness of $\Omega$ and closedness of ${\mathcal S}_{0}$ one has $d(\phi({\mathcal S}_{0}), \partial \Omega)>0$ for $\phi \in {\mathcal C}$. \par
\ \par
Now the tangent space of ${\mathcal C}$ is described by the following proposition.
\begin{Proposition} \label{Prop:UnIIsatisfaitcaC}
Let $(\tau,\eta) \in {\mathcal C}$. Then using the notations 
\begin{equation} \label{NotationsC}
x_B := \tau (x_0), \ {\mathcal S}  := \tau  ({\mathcal S}_{0}) \text{ and } \ {\mathcal F}  := \Omega \setminus {\mathcal S}, 
\end{equation}
we have
\begin{align}
\nonumber
T_{(\tau,\eta)} {\mathcal C} =  \bigg \{ (\sigma,\mu) \in T_{\tau}SE(3) &\times C^{1,\lambda}({\mathcal F}_{0};\R^{3}) \ \Big/ \\ %\\
\label{UDiv}
& \div(U_{\eta}[\mu])=0 \ \text{ in }\ {\mathcal F}, \\
\label{UAuBord}
& U_{\eta}[\mu](x) \cdot n(x) = (L_{\tau}[\sigma]  + R_{\tau}[\sigma] \times (x -x_B)) \cdot n(x) \ \text{ on }\ \partial {\mathcal S} , \\
\label{UAuBordFixe}
& U_{\eta}[\mu](x) \cdot n(x) =0 \ \text{ on }\ \partial \Omega \bigg\}.
\end{align}
As before, $n$ denotes the normal unit vector on $\partial \Omega$ and $\partial {\mathcal S}$, pointing outside the fluid.
\end{Proposition}
The proof of Proposition \ref{Prop:UnIIsatisfaitcaC} will be given in the next section. We will represent $(\sigma,\mu)$ by ${\mathcal U} \circ \phi: \Omega \rightarrow \R^{3}$ where ${\mathcal U}:\Omega \rightarrow \R^{3}$ is given by ${\mathcal U}_{|{\mathcal S}} \circ \tau= \sigma$ and ${\mathcal U}_{|{\mathcal F}}= U_{\eta}[\mu]$. In the same way as $\phi$, ${\mathcal U}$ can be discontinuous. \par
\ \par
The manifold ${\mathcal C}$ can be endowed by the following Riemannian metric:  for any  $\phi \in {\mathcal C}$, for any ${\mathcal U}_1 \circ \phi$, ${\mathcal U}_2 \circ \phi \in T_{\phi  } {\mathcal C}$,
\begin{align} \nonumber
\langle {\mathcal U}_1 \circ \phi , {\mathcal U}_2 \circ \phi  \rangle_{\phi }:=  
  \int_{  \Omega } \rho_0 ( \phi^{-1}) \, {\mathcal U}_1  \cdot {\mathcal U}_2 \, dx  .
\end{align}
Above we used the density $\rho_0 $ defined on $\Omega$ by:
\begin{equation*}
\rho_0 (x) = \left\{ \begin{array}{l}
 1 \text{ in } {\mathcal F}_{0}, \\
 \rho_{{\mathcal S}_0}(x) \text{ in } {\mathcal S}_{0}.
\end{array} \right.
\end{equation*}
Splitting the fluid flow and the solid one this reads: for any $(\sigma_1 ,\mu_1) , (\sigma_2 ,\mu_2 )$ in $T_{(\tau,\eta)} {\mathcal C}$, we have
\begin{align} \label{metrique}
\langle (\sigma_1 ,\mu_1) , (\sigma_2 ,\mu_2 ) \rangle_{(\tau,\eta)} =  m L_{\tau}[\sigma_{1}]  \cdot L_{\tau}[\sigma_{2}]  +  {\mathcal J}[\tau] R_{\tau}[\sigma_{1}]  \cdot  R_{\tau}[\sigma_{2}] +   \int_{ \eta ({\mathcal F}_{0}) } U_{\eta}[\mu_{1}] \cdot  U_{\eta}[\mu_{2}]  \, dx ,
\end{align}
where ${\mathcal J}[\tau]$ is the inertia matrix deduced from the initial one ${\mathcal J}_0$ by the rigid transformation $\tau$ that is 
\begin{equation} \label{Jtau}
\mathcal{J} [\tau] =Q[\tau] \mathcal{J}_0 Q[\tau]^{*},
\end{equation}
where $Q[\tau]$ is the rotation matrix canonically associated to $\tau$, that is, its linear part. \par
\ \par
Let us stress that this metric defines a weaker topology than the original one on ${\mathcal C}$. \par
% Let us also observe that the group $SE(3) \times \text{Sdiff}^{+}  ({\mathcal F}_{0} )$ acts from the right on the manifold ${\mathcal C}$ by composition and leaves its metric invariant, which corresponds to the particle relabelling symmetry.
%
%
\ \\
{\bf Curves of configurations.} We now turn to time-dependent displacements. Let us be given $(\tau_{0},\eta_{0})$ and $(\tau_{1},\eta_{1})$ in ${\mathcal C}$. We introduce
\begin{align} \nonumber
{\mathcal L} :=  \Big\{ \ (\tau,\eta) &\in C^{1}([0,T];{\mathcal C}), \text{ such that:}  \\
\nonumber
&{i.} \ \ \tau(0)=\tau_{0}, \ \ \eta(0)=\eta_{0}, \\
&{ii.} \ \ \tau(T)=\tau_{1}, \ \ \eta(T)=\eta_{1}
\ \Big\}.
\label{defL}
\end{align}
It is easy to verify that ${\mathcal L}$ is a submanifold of a manifold modelled on the  Banach space
$E_{T}:= C^{1}([0,T]; \mathfrak{se}(3)) \times C^{1,\lambda}([0,T] \times {\mathcal F}_{0}; \R^{3})$. \par
The tangent space of this manifold is described by the following proposition.
\begin{Proposition} \label{Prop:UnIIsatisfaitcaG}
Let $(\tau,\eta) \in {\mathcal L}$. Then using the notations 
\begin{equation} \label{NotationsL}
x_B (t) := \tau_{t} (x_0), \ {\mathcal S} (t)  := \tau_{t}({\mathcal S}_{0}) \ \text{ and } \ {\mathcal F}(t):=\Omega \setminus {\mathcal S}  (t),
\end{equation}
we have
\begin{align}
\nonumber
T_{(\tau,\eta)} {\mathcal L} =  \Big\{ &(\sigma,\mu) \in C^{1}([0,T]; T_{\tau_{t}}SE(3)) \times C^{1,\lambda}([0,T] \times {\mathcal F}_{0}; \R^{3}) \ \Big/  \\
\label{Nulsen0etT}
& \sigma(0)=\sigma(T)=0 \text{ and } \mu(0)=\mu(T)=0, \\
%&U_{\eta}[\mu](0)=U_{\eta}[\mu](T)= 0 \text{ and } L_{\tau}[\sigma](0)= R_{\tau}[\sigma](0)=L_{\tau}[\sigma](T)= R_{\tau}[\sigma](T)=0, \\
%
\label{UDiv2}
& \div(U_{\eta}[\mu])=0 \text{ in } {\mathcal F}(t) \text{ for each } t \in [0,T], \\
\label{UAuBord2}
& U_{\eta}[\mu] \cdot n =(L_{\tau}[\sigma] + R_{\tau}[\sigma] \times [x-x_B (t)])  \cdot  n \text{ on } \partial {\mathcal S}(t)
\ \text{ for each } t \in [0,T], \\
\label{UAuBordFixe2}
& U_{\eta}[\mu] \cdot  n =0 \text{ on } \partial \Omega \text{ for each } t \in [0,T] \Big\}.
\end{align}
\end{Proposition}
\begin{Remark}
Here we make the abuse of notations $C^{1}([0,T];T_{\tau_{t}}SE(3))$, since $T_{\tau_{t}}SE(3)$ actually depends on $t$. We consider $\sigma$ as the section of a vector bundle rather than as a function. One can for instance interpret $C^{1}([0,T];T_{\tau_{t}}SE(3))$ as the set of mappings $\sigma$ such that $(t,x) \mapsto \sigma(t,\tau_{t}^{-1}(x))$ is in $C^{1}([0,T];\mathfrak{se}(3))$. Also, we sometimes drop the dependence of the objects on $t$ in order to simplify the notations.
\end{Remark}
\ \par
\noindent
Proposition \ref{Prop:UnIIsatisfaitcaG} is the time-dependent counterpart of  Proposition \ref{Prop:UnIIsatisfaitcaC}. We do not provide a proof since it is only a matter of adapting the proof of Proposition \ref{Prop:UnIIsatisfaitcaC} with a harmless parameter. 
The only new point is to observe that the extremities of the curves being prescribed (the conditions i) and ii) in \eqref{defL}) the fields $\sigma$ and $\mu$ vanish when $t=0$ or $T$. \par
\ \par
As previously we will represent $(\tau,\eta) \in {\mathcal L}$  by $\phi : [0,T] \times \Omega \rightarrow \Omega$ such that for any $t \in [0,T]$, 
\begin{equation*}
\phi (t,\cdot) |_{{\mathcal S}_{0}}= \tau  (t,\cdot) \ \text{ and } \ \phi (t,\cdot) |_{{\mathcal F}_{0}}=\eta  (t,\cdot).
\end{equation*}
We represent $(\sigma,\mu) \in T_{(\tau,\eta)} {\mathcal L}$ by ${\mathcal U} \circ \phi: [0,T] \times  \Omega \rightarrow \R^{3}$ where ${\mathcal U}: [0,T] \times \Omega \rightarrow \R^{3}$ is given, for any $t \in [0,T]$, by 
\begin{equation*}
{\mathcal U}(t,\cdot) |_{{\mathcal S}(t)} \circ \tau (t,\cdot) = \sigma (t,\cdot) \ 
\text{ and } \ {\mathcal U}(t,\cdot) |_{{\mathcal F}(t)}= U_{\eta}[\mu] (t,\cdot) . 
\end{equation*}
\subsubsection{The geodesic interpretation of the  motion of a rigid body immersed in an incompressible perfect  fluid} 
\label{GeodesicInterpretation}
Here we consider geodesics as critical points of the following action on the manifold ${\mathcal L}$:
\begin{equation} \label{DefA}
{\mathcal A}(\phi) := \frac{1}{2} \int_{[0,T]  \times \Omega }  \rho_{{0}}  |\partial_t  \phi |^2 \, dx\, dt. 
\end{equation}
%
%\\  &= \frac{1}{2} \int_{[0,T]} \Big( m |\ell(t)|^2  + \langle J[\tau(t)] r(t), r(t) \rangle +  \int_{{\mathcal F}_0 }  |\partial_t \eta|^2 dx  \Big) \, dt,
%
We see that the action is obtained by integrating the squared norm (associated to the metric of  ${\mathcal C}$) of the tangent vector to the curve $\phi$, that is
\begin{equation*}
{\mathcal A}(\phi) = \frac{1}{2} \int_{[0,T] } \langle  \partial_t  \phi (t,\cdot),\partial_t  \phi  (t,\cdot) \rangle_{\phi (t,\cdot)} \, dt .
\end{equation*}
Separating the fluid and the body parts in the integral, and using \eqref{meti}, we see that
\begin{equation*}
{\mathcal A}(\phi)= \frac{1}{2} \int_{[0,T]} \Big( m |\ell(t)|^2  + {\mathcal J}[\tau(t)] r(t) \cdot r(t) +  \int_{\eta_{t}({\mathcal F}_{0}) }  |u(t,x)|^2 \, dx  \Big) \, dt,
\end{equation*}
where
\begin{equation} \label{LRU}
\ell:=L_{\tau}[\partial_{t} \tau], \ r:=R_{\tau}[\partial_{t} \tau] \ \text{ and } \ u:=U_{\eta}[\partial_{t} \eta].
\end{equation}
In this writing, we recognize the integral over time of the kinetic energy of the fluid-body system. \par
Moreover, going back to \eqref{DefA}, since the action is a continuous quadratic form on ${\mathcal L}$, we deduce that ${\mathcal A}$ is differentiable on ${\mathcal L}$ with
\begin{align} \label{diff}
 D {\mathcal A} (\phi) \cdot ( {\mathcal U} \circ \phi )  = 
\int_{[0,T]  \times \Omega }  \rho_{{0}} \,  \partial_t  \phi  \cdot  \partial_t   ( {\mathcal U} \circ \phi )  \, dx\, dt 
= \int_{[0,T]} \langle  \partial_t  \phi (t,\cdot),\partial_t ({\mathcal U} \circ \phi)  (t,\cdot) \rangle_{\phi (t,\cdot)} \, dt .
\end{align}
This leads us to the following natural definition.
\begin{Definition}
\label{geocca}
We say that $\phi \in {\mathcal L}$ is a geodesic on ${\mathcal L}$ if for any  ${\mathcal U} \circ \phi  \in T_{\phi} {\mathcal L}$, $ D {\mathcal A} (\phi) \cdot ({\mathcal U} \circ \phi ) = 0$.
\end{Definition}
\subsection{Equivalence of the two points of view} 
\label{Equi}
The main result of this paper is the following.
\begin{Theorem}
\label{Theo:Equi}
If $(u,x_B , r) $ is a  classical solution of the PDEs formulation on $[{0}, T]$ then $(\tau,\eta) $ defined by formulas \eqref{flowS} and \eqref{flowF} is a  geodesic on ${\mathcal L}$. \par
Conversely, let $(\tau,\eta) \in {\mathcal L}$ be a geodesic. Then   $(u,x_B , r) $ where $(\ell,r,u)$ is defined by \eqref{LRU} 
and $x_{B}$ is obtained by \eqref{CMci} and \eqref{DefEll}, is a classical solution of the PDEs formulation on $[{0}, T]$ (in the sense of Definition \ref{germi}).
\end{Theorem}
Theorem \ref{Theo:Equi} will be proved in Subsection \ref{EquiProof}. \par
\ \par
\begin{Remark}
Let us define the length of a curve $\phi \in {\mathcal L}$:
\begin{equation*}
\Lambda(\phi) := \int_{[0,T]}  ( \langle  \partial_t  \phi (t,\cdot),\partial_t  \phi  (t,\cdot) \rangle_{\phi (t,\cdot)})^\frac{1}{2} \, dt ,
\end{equation*}
and consider
\begin{align} \nonumber
d &:=  \inf \Lambda(\phi) \\
&= \inf \int_{[0,T]} \left( |\ell(t)|^2  + {\mathcal J}(t) r(t) \cdot  r(t)  +  \int_{\eta_{t}({\mathcal F}_0) }  |u(t,x)|^2 \, dx \right)^{1/2} \, dt ,
\end{align}
where the infimum is performed over  $\phi = (\tau,\eta) \in {\mathcal L}$. We should say that $d$ is the geodesic distance between the configurations $(\tau_{0}, \eta_{0})$ and $(\tau_{1},\eta_{1})$ of  ${\mathcal C}$. 
If $(\tau,\eta) \in {\mathcal L}$ realizes this infimum and is parametrized by $t$ in such a way that the energy does not depend on time then  $(\tau,\eta)$ also minimizes the action ${\mathcal A}$ over ${\mathcal L}$.
Conversely, by the conservation of energy, any geodesic is parameterized proportionnaly to arc length.
\end{Remark}

Let us mention here two open problems. \par
\ \\
{\bf Open Problem 1.} Is it possible to prove that for $T$ small enough and $(\tau,\eta) \in {\mathcal L}$ such that the associated  $(u,x_B , r) $ is a  classical solution of the PDEs formulation on $[{0}, T]$, one has for any $({\tilde \tau},{\tilde \eta}) \in {\mathcal L}$, ${\mathcal A} (\tau,\eta) \leqslant {\mathcal A} ({\tilde \tau},{\tilde \eta})$, with equality if and only if  $({\tilde \tau},{\tilde \eta}) = (\tau,\eta)$? This should extend the result obtained by  Brenier  cf. \cite{Brenier} in the case of a fluid without body. \par
%
% En temps long, je crois que Schnirelman dit no way.
\ \\
{\bf Open Problem 2.} Is it possible to adapt the strategy that Ebin and Marsden used in \cite{EbinMarsden} in the case of a fluid alone to the case with a body, that is to prove the existence of a free torsion connexion and of some geodesics by using parallel transport, and then to prove that these solutions also solve the PDE formulation? \par
\ \\
Let us also mention studies connected with the stability properties of the system. In the case of a fluid alone, Arnold \cite{Arnold} uses a notion of Riemannian curvature to investigate the stability of two-dimensional stationary flows. In the case considered here of a fluid-body system, this stability was studied by Ilin and Vladimirov \cite{VladimirovIlin2D,VladimirovIlin3D}.
%
%
%%%%%%%%%%%%%%%%%%%%%%%%%%%%%%%%%%%%%%%%%%%%%%%%%%%%%%
%
%
\section{Proofs}
\label{Sec:Proofs}
In this section, we prove the claims of Section \ref{Intro}. The existence of a solution of the PDE system will not be needed.
\subsection{Proof of Proposition \ref{Prop:UnIIsatisfaitcaC}}
\label{TS}
{\bf 1.} Let us first show that any $(\sigma,\mu)$ satisfying \eqref{UDiv}-\eqref{UAuBordFixe} belongs to $T_{(\tau,\eta)} {\mathcal C}$.  First,  we define $\theta_{\mathcal S} \in C^{1}((-\varepsilon,\varepsilon);SE(3))$ by
\begin{equation*}
\partial_s \theta_{\mathcal S} (s) [x]= v_{\mathcal S} ( s, \theta_{\mathcal S} (s) [x] ) , \quad \theta_{\mathcal S} (0) [x]=\tau [x] ,
\end{equation*}
where 
\begin{equation*}
v_{\mathcal S} (s,x) := L_{\tau}[\sigma] + R_{\tau}[\sigma] \times (x-x_B (s)) , \quad x_B (s) := x_B + s L_{\tau}[\sigma] ,
\end{equation*}
with $\varepsilon>0$ small; in particular we can require that for all $s \in (-\varepsilon,\varepsilon)$, one has
\begin{equation*}
{\mathcal S}_{s}:=\theta_{\mathcal S}(s)[{\mathcal S}] \subset \Omega.
\end{equation*}
\begin{Lemma}
There exists a family 
$\psi_{s} : {\mathcal F} \rightarrow {\mathcal F}_{s}:= \Omega \setminus {\mathcal S}_{s}$, for $s \in (-\varepsilon,\varepsilon)$, smooth in its arguments, such that 
\begin{gather}
\label{F1}
 \psi_{0}=\Id_{{\mathcal F}}, \\
\label{F2}
\forall s \in (-\varepsilon,\varepsilon),  \ \ \psi_{s}= \theta_{\mathcal S}(s) \text{ in a neighborhood of } {\mathcal S}, \\
\label{F3}
\forall s \in (-\varepsilon,\varepsilon), \ \ \psi_{s}= \Id \text{ in a neighborhood of } \partial \Omega, \\
\label{F4}
\forall s \in (-\varepsilon,\varepsilon) , \ \forall x \in {\mathcal F}, \ \ \det(\nabla_{x} \, \psi_{s})=1.
\end{gather}
\end{Lemma}
\begin{proof}
We define 
\begin{equation} \label{DefV}
V(s,x) := -\frac{1}{2} \curl \big( \phi(x) (x\times L_{\tau}[\sigma] + \| x - x_B (s) \|^2  R_{\tau}[\sigma] ) \big),
\end{equation}
where $\phi$ is a smooth function equal to $1$ in a neighborhood of ${\mathcal S}$ and $0$  in a neighborhood of $\partial \Omega$. We associate $\psi_{s}$ as the solution of 
\begin{equation*}
\partial_s \psi_{s} (x)= V( s, \psi_{s} (x)) , \quad   \psi_{0} (x)=x .
\end{equation*}
Equations \eqref{F1} and \eqref{F3} are straightforward. Equation \eqref{F2} follows from the fact that $V$ coincides with $v_{\mathcal S}$ in some neighborhood of $\partial {\mathcal S}$, reducing $\varepsilon$ if necessary. Finally, since $V$ is clearly divergence-free, \eqref{F4} follows from Liouville's theorem.
\end{proof}
%
%
%\ \par
%
%
Now we construct the vector field:
\begin{equation} \label{defvt}
v_{\mathcal F} (s, \cdot) := V(s,\cdot ) + (\psi_{s})_{*}( U_{\eta}[\mu] - V( 0 ,\cdot ) ) \ \text{ in } {\mathcal F}_{s}.
\end{equation}
(Recall that the push-forward is defined as $(\psi_{*} v  )(y) =  d\psi_{\psi^{-1}(y)} [v (\psi^{-1}(y))] $.) \par
Now we introduce a corresponding flow $\theta_{\mathcal F}$ according to the $s$-variable, starting from $\eta$:
\begin{equation*}
\left\{ \begin{array}{l}
 \partial_{s} \theta_{\mathcal F} (s,x) = v_{\mathcal F} (s,\theta_{\mathcal F}(s,x)), \ \text{ in } (-\varepsilon,\varepsilon) \times {\mathcal F}_{0}, \\
 \theta_{\mathcal F}(0,x) = \eta(x) \ \text{ in } {\mathcal F}_{0}.
\end{array} \right.
\end{equation*}
(One can for instance smoothly extend $v_{{\mathcal F}}(s,\cdot)$ to $\R^{n}$ to define this flow in a standard way.) \par
\ \par
It remains to check that $\Theta:= (\theta_{\mathcal S},\theta_{\mathcal F}) \in C^{1}((-\varepsilon,\varepsilon);{\mathcal C})$, that $\Theta(0)=(\tau,\eta)$ and that $\Theta '(0)=  (\sigma , \mu)$. The two latter claims come directly from the construction.\par
Let us prove that $\theta_{\mathcal F}$ is volume-preserving. To that purpose, we first notice that $\div(U_{\mu})=\div(V(s))=0$. Using the fact that the push-forward of a divergence-free vector field by a diffeomorphism with unit Jacobian determinant is still divergence-free (see for instance \cite[Proposition 2.4]{InoueWakimoto}), \eqref{F4} and \eqref{defvt}, we deduce that $ v_{\mathcal F} (s)$ is divergence-free. Hence it follows that $\theta_{\mathcal F}$ is volume-preserving by Liouville's theorem. \par
The main point, that is that, for any $s \in (-\varepsilon,\varepsilon)$, $\theta_{\mathcal F}(s)$ sends ${\mathcal F}_{0}$ to ${\mathcal F}_{s}$ can be seen as follows. It suffices to prove that 
\begin{gather}
\label{Tangence1}
v_{\mathcal F} (s,\cdot).n = v_{\mathcal S} (s,\cdot).n \ \text{ on } \partial {\mathcal S}_{s}, \\
\label{Tangence2}
v_{\mathcal F} (s,\cdot).n = 0 \ \text{ on } \partial \Omega.
\end{gather}
Using \eqref{F3} and \eqref{DefV}, we see that $v_{\mathcal F} = U_{\eta}[\mu]$ in some neighborhood of $\partial \Omega$, so \eqref{Tangence2} is a consequence of \eqref{UAuBordFixe}. For what concerns \eqref{Tangence1}, from \eqref{UAuBord} we see that $U_{\eta}[\mu] - v_{\mathcal S} (0,\cdot)$ is tangent to $\partial {\mathcal S}$. It follows that $(\psi_{s})_{*}(U_{\eta}[\mu] - v_{\mathcal S} (0,\cdot))$ is tangent to $\partial {\mathcal S}_{s}$, which gives \eqref{Tangence1}. \par
The other requirements, namely that $\theta_{\mathcal F}$ is orientation preserving and has the claimed regularity, are clearly satisfied. \par
\ \par
\noindent
{\bf 2.} Reciprocally, given $(\sigma,\mu)$ in $T_{(\tau,\eta)} {\mathcal C}$, by definition of the tangent space, there is $\Theta= (\theta_{{\mathcal S}},\theta_{{\mathcal F}}) \in C^{1}((-\varepsilon,\varepsilon); {\mathcal C})$  such that $\Theta(0)=(\tau,\eta)$ and that $\Theta '(0)=  (\sigma , \mu)$. Let us check that 
\eqref{UDiv}-\eqref{UAuBordFixe} are satisfied. \par
First, it is obvious that $\sigma=\theta'_{{\mathcal S}}(0) \in T_{\tau} SE(3)$. Next one can define
\begin{equation*}
v_{{\mathcal F}}(s,x):= \partial_{s} \theta_{{\mathcal F}}(s, \theta^{-1}_{{\mathcal F}}(s,x)),
\end{equation*}
as a vector field on ${\mathcal F}_{s}$. Then as
\begin{equation} \label{TF}
\partial_{s} \theta_{\mathcal F}(s, x) = v_{{\mathcal F}}(s,\theta_{\mathcal F}(s,x)),
\end{equation}
and since $\theta_{{\mathcal F}}$ is volume-preserving, by Liouville's theorem we infer that $\div v_{{\mathcal F}}=0$. Using again \eqref{TF} and the fact that $\theta_{{\mathcal F}}$ sends ${\mathcal F}_{0}$ to ${\mathcal F}_{s}$, we also easily infer \eqref{UAuBord} and \eqref{UAuBordFixe}. \eop
\subsection{Proof of Theorem \ref{Theo:Equi}}
\label{EquiProof}

We start with the following lemma.
%Here again one can split the fluid flow and the solid one to obtain, with the previous notations, the following.
%
\begin{Lemma}
\label{plusimple}
For any $\phi =(\tau,\eta) \in {\mathcal L}$, for any ${\mathcal U} \circ \phi \in T_{\phi} {\mathcal L}$,
\begin{align} \nonumber
&D {\mathcal A} (\phi) \cdot ( {\mathcal U} \circ \phi )  =  \\
\nonumber
&\quad \int_{[0,T]} \Big( m \ell (t) \cdot  \partial_{t} L_{\tau}[\sigma] (t)
+ {\mathcal J} [\tau] r (t) \cdot  \partial_{t} R_{\tau}[\sigma]  (t)
+ \int_{{\mathcal F}_0 }  \partial_t \eta (t,x)  \cdot \partial_t ( {\mathcal U}(t,\eta(t,x)) ) \, dx \Big) \, dt,
\end{align}
where $\ell$ and $r$ are given by \eqref{LRU}, $Q[\tau]$ is the linear part of $\tau$ and ${\mathcal J
}[\tau]$ is given by \eqref{Jtau}.
\end{Lemma}
\begin{proof}
Let us split the fluid flow and the solid one in \eqref{diff} and prove that 
\begin{align}
\label{idplusimple}
 \int_{[0,T]  \times  {\mathcal S}_0 }  \rho_{{\mathcal S}_0} \, \partial_t  \tau  \cdot  \partial_t  \sigma   \, dx\, dt  =
\int_{[0,T]} \Big( m \ell  \cdot  \partial_{t} L_{\tau}[\sigma]  
 +   {\mathcal J} [\tau] r  \cdot  \partial_{t} R_{\tau}[\sigma]   \Big) .
\end{align}
Now in order to prove \eqref{idplusimple} we first perform a change of variable to get
%observe that the contribution of the solid flow is
%
\begin{align} \label{maissi}
 \int_{[0,T]  \times  {\mathcal S}_0 }  \rho_{{\mathcal S}_0} \, \partial_t  \tau  \cdot  \partial_t  \sigma   \, dx\, dt  =
 \int_{[0,T] }  \int_{{\mathcal S}(t) }  \rho_{{\mathcal S}_0} (\tau_t^{-1}(x))  \, \partial_t  \tau (t,\tau_t^{-1}(x))   \cdot  \partial_t  \sigma(t,\tau_t^{-1}(x))   \, dx\, dt .
\end{align}
%
%Above the notations $\tau_t^{-1}$ and $\eta_t^{-1}$ stand for the respective inverse of the functions $\tau (t,\cdot)$ and  $\eta (t,\cdot)$. 
Now, by definition of $\ell$ and $r$ we have 
\begin{align} \label{maissiphi}
 \partial_t  \tau(t, \tau_t^{-1}(x))  = \ell + r \times (x-x_B ) .
 \end{align}
In particular, taking the linear part of these affine transformations, and using that the linear part of $\partial_{t} \tau$ is $\partial_{t} Q[\tau]$, this entails that 
\begin{align} \label{maissiQ}
  \partial_t Q [\tau]=  r \times Q[\tau] .
 \end{align}
On the other hand, by definition of $L_{\tau}[\sigma]$ and $R_{\tau}[\sigma]$, we have 
\begin{align} \nonumber
\sigma  = L_{\tau}[\sigma] + R_{\tau}[\sigma] \times Q[\tau] (x-x_0 ) ,
\end{align}
so that, differentiating in time and using \eqref{maissiQ}, we obtain
\begin{align} \nonumber
\partial_t  \sigma  = \partial_t  L_{\tau}[\sigma] + (\partial_t  R_{\tau}[\sigma]) \times Q[\tau] (x-x_0 ) + 
R_{\tau}[\sigma] \times (r  \times Q[\tau] (x-x_0 )) .
\end{align}
Operating $\tau_t^{-1}$ on the right for both sides of the previous equality, we get 
\begin{align} \label{maissitau}
\partial_t  \sigma(t,\tau_t^{-1}(x))  = \partial_{t} L_{\tau}[\sigma]  +  \partial_{t} R_{\tau}[\sigma]  \times (x-x_B ) + R_{\tau}[\sigma]  \times \Big( r \times (x-x_B ) \Big) .
\end{align}
We now plug \eqref{maissiphi} and \eqref{maissitau} into the right hand side of 
\eqref{maissi} to obtain
\begin{align}
\label{moinslong}
 \int_{[0,T]  \times  {\mathcal S}_0 }  \rho_{{\mathcal S}_0} \partial_t  \tau  \cdot  \partial_t  \sigma   \, dx\, dt  =
\int_{[0,T]}  ( I_1 (t) + I_2 (t) ) dt ,
\end{align}
with
\begin{align*}
I_1 (t) &:=  \int_{{\mathcal S}(t) }  \rho_{{\mathcal S}_0} (\tau_t^{-1}(x))
\Big(  \ell (t) + r(t) \times (x-x_B (t) )  \Big)  \cdot  
\Big(  \partial_{t} L_{\tau}[\sigma] (t) +  \partial_{t} R_{\tau}[\sigma] (t) \times (x-x_B (t))  \Big)
\, dx , \\  
I_2 (t) &:=  \int_{{\mathcal S}(t) }  \rho_{{\mathcal S}_0} (\tau_t^{-1}(x))
\Big(  \ell(t) + r(t) \times (x-x_B(t) )  \Big)    \cdot
\Big[  R_{\tau}[\sigma] (t) \times \Big( r(t) \times (x-x_B (t)) \Big) \Big] \, dx .
\end{align*}
We use the identity \eqref{meti} with 
\begin{equation*}
( \ell_1, \, r_1, \, \ell_2, \, r_2 ) = (\ell (t),\,  r (t), \, \partial_{t} L_{\tau}[\sigma] (t), \, \partial_{t} R_{\tau}[\sigma] (t) ),
\end{equation*}
to get 
\begin{align}
\label{hihun} 
 I_1 (t) = m \ell (t) \cdot  \partial_{t} L_{\tau}[\sigma] (t) 
 +   {\mathcal J} [\tau(t)] r (t) \cdot  \partial_{t} R_{\tau}[\sigma]  (t) .
 \end{align}
Finally we observe that 
\begin{align}
\nonumber
 I_2 (t) &= \int_{{\mathcal S}(t) }  \rho_{{\mathcal S}_0} (\tau_t^{-1}(x)) \, \ell (t) \cdot
\Big(  R_{\tau}[\sigma] (t) \times \big[ r (t) \times (x-x_B (t)) \big] \Big) \, dx  \\ 
\nonumber
&= \ell (t)   \cdot    \bigg(  R_{\tau}[\sigma] (t) \times
\Big( r(t) \times \Big[ \int_{{\mathcal S}(t) } \rho_{{\mathcal S}_0} (\tau_t^{-1}(x))  (x-x_B (t))  \,dx \Big] \Big)
 \bigg) \\
& = 0 , \label{hideux} 
  \end{align}
according to \eqref{premeti}. Combining \eqref{moinslong}, \eqref{hihun} and \eqref{hideux} yields \eqref{idplusimple}.
\end{proof}
\ \par
\noindent
\begin{proof}[Proof of Theorem \ref{Theo:Equi}] \ \par
\ \par
\noindent
{\bf 1.} Let $(u,x_B , r) $ is a  classical solution of the PDEs formulation on $[{0}, T]$ and let $(\tau,\eta) $ their respective flows  given by formulas \eqref{flowS} and \eqref{flowF}. Equation  (\ref{Euler1a2}) reads 
\begin{equation} \label{Leonard}
\displaystyle \frac{\partial^2 \eta}{\partial t^2} + (\nabla p ) \circ \eta = 0 ,  \ \text{for}  \ x\in   \mathcal{F}_0 .
\end{equation}
Let us now consider $(\sigma,\mu)$ in $T_{(\tau,\eta)} {\mathcal L}$. Using \eqref{Leonard} and performing a change of variables, we get
\begin{equation*}
- \int_{{\mathcal F}_{0}  }   \frac{\partial^2 \eta}{\partial t^2} \cdot U_{\eta}[\mu] \circ \eta  \, dx =
- \int_{{\mathcal F} (t)  }  \nabla p \cdot U_{\eta}[\mu]   \, dx.
\end{equation*}
Now using \eqref{UDiv2}, \eqref{UAuBordFixe2} and Green's formula we deduce
\begin{equation*}
- \int_{{\mathcal F} (t)  }   \frac{\partial^2 \eta}{\partial t^2} \cdot U_{\eta}[\mu] \circ \eta  \, dx
= \int_{\partial \mathcal{S}(t)} p \, U_{\eta}[\mu]  \cdot n  \, d \Gamma.
\end{equation*}

Now using (\ref{UAuBord2}) and then (\ref{Solide1})-(\ref{Solide2}), we obtain
\begin{eqnarray*}
- \int_{{\mathcal F} (t)  }   \frac{\partial^2 \eta}{\partial t^2} \cdot U_{\eta}[\mu] \circ \eta  \, dx
&= & \int_{\partial \mathcal{S}(t)} p \, \big ( L_{\tau}[\sigma] + R_{\tau}[\sigma] \times (x-x_{B}(t)) \big)\cdot n  \, d \Gamma(x) \\  
&= & m x''_B (t) \cdot L_{\tau}[\sigma]   + ({\mathcal J} [\tau ] (t) r (t))'  \cdot R_{\tau}[\sigma].
\end{eqnarray*}
Then we integrate by parts over $[0,T]$ and conclude with Lemma \ref{plusimple} that $(\tau,\eta) $ is a  geodesic on ${\mathcal L}$. \par
\ \par
\noindent
{\bf 2.} Conversely, let $\phi \ \in {\mathcal L}$ be a geodesic. Using again Lemma \ref{plusimple}, it means that for any $(\sigma,\mu) \in T_{(\tau,\eta)} {\mathcal L}$, one has
\begin{equation} \label{Var}
\int_{[0,T]} \Big( m \ell(t) \cdot  L_{\tau}[\sigma] (t)   + {\mathcal J} [\tau(t) ] r (t) \cdot R_{\tau}[\sigma]  (t) 
+  \int_{{\mathcal F}_0 }  \partial_t \eta (t,x) \cdot  \partial_t ( U_{\tau}[\mu](t,\eta(t,x))  )  \, dx \Big) \, dt =0.
\end{equation}
We first use \eqref{Var} with $\sigma=0$. We consider  $w : [0,T] \times {\mathcal F}_{t} \rightarrow \R^{3}$ satisfying $\div(w)=0$ in ${\mathcal F}_{t}$, $w \cdot n=0$ on $\partial {\mathcal F}_{t}$ and $w(0,\cdot)=0$ and $w(T,\cdot)=0$. Then $(\sigma,\mu) \in T_{(\tau,\eta)} {\mathcal L}$ for $\sigma=0$ and $\mu(t,x):= w(t,\eta(t,x))$. Consequently, one has
\begin{align*}
0 &= \int_{[0,T]} \int_{{\mathcal F}_0 }  \partial_t \eta (t,x) \cdot  \partial_t (w(t,\eta(t,x))) \, dx  \, dt \\
&= - \int_{[0,T]} \int_{{\mathcal F}_0 }  \partial^{2}_t \eta (t,x) \cdot  w(t,\eta(t,x)) \, dx  \, dt \\
&= - \int_{[0,T]} \int_{{\mathcal F}_t }  \partial^{2}_t \eta (t,\eta^{-1}(t,x)) \cdot  w(t,x) \, dx  \, dt.
\end{align*} 
It follows that $\partial^{2}_t \eta (t,\eta^{-1}(t,\cdot))$ is a gradient field in ${\mathcal F}_{t}$, that we denote $-\nabla p$. Now going back to \eqref{Var} we get that for general $(\sigma,\mu) \in T_{(\tau,\eta)} {\mathcal L}$, one has
\begin{equation*}
\int_{[0,T]} \Big( m \ell(t) \cdot  \partial_{t} L_{\tau}[\sigma] (t)  +  {\mathcal J} [\tau(t) ] r (t) \cdot  \partial_{t} (R_{\tau}[\sigma]) (t)
  +  \int_{{\mathcal F}_t }  \nabla p(t,x)  \cdot  U_{\eta}[\mu](t,x) \, dx \Big) \, dt =0.
\end{equation*}
Using Green's formula for the last integral and \eqref{UAuBord2}-\eqref{UAuBordFixe2}, we deduce that for any $\sigma$,
\begin{multline*}
\int_{[0,T]} \bigg( m   \ell(t) \cdot  \partial_{t} L_{\tau}[\sigma] (t)
+  {\mathcal J} [\tau(t) ] r (t)  \cdot  \partial_{t} (R_{\tau}[\sigma]) (t)  \\
+ L_{\tau}[\sigma]  \cdot  \int_{\partial {\mathcal F}_t } p(t,x) n(x)  \, d \Gamma 
+ R_{\tau}[\sigma]  \cdot  \int_{\partial {\mathcal F}_t } p(t,x) (x-x_{B}(t)) \times n(x)  \, d \Gamma \bigg) \, dt =0.
\end{multline*}
Now we integrate by parts in time the first two terms. Since this is valid for any $\sigma$, hence for any $L_{\tau}[\sigma]$ and $R_{\tau}[\sigma]$, we infer \eqref{Solide1}-\eqref{Solide2}. Equations \eqref{Euler3a2}-\eqref{Euler3b} then follow from the very definition of ${\mathcal L}$.
\end{proof}
\ \par
\noindent
{\bf Acknowledgements.} The authors were partially supported by the Agence Nationale de la Recherche, Project CISIFS, grant ANR-09-BLAN-0213-02.

\end{document}